# A Scheduling Problem for Hospital Operating Theatre

Suliadi F. Sufahani[1], Siti N. A. Mohd Razali[2], Zuhaimy Ismail[3]

[1,2]Department Of Mathematic and Sciences, Faculty of Science, Technology and Human Development, University Tun Hussein Onn Malaysia, 86400 Parit Raja, Batu Pahat, Johor, Malaysia
suliadi@uthm.edu.my     asyikinr@uthm.edu.my

[3]Deparmtent Of Mathematic, Faculty of Science, University Teknologi Malaysia, Skudai, Johor Bahru, Johor, Malaysia
zuhaimyi@yahoo.com

April 2012

This paper provides a classification of real scheduling problems. Various ways have been examined and described on the problem. Scheduling problem faces a tremendous challenges and difficulties in order to meet the preferences of the consumer. Dealing with scheduling problem is complicated, inefficient and time-consuming. This study aims to develop a mathematical model for scheduling the operating theatre during peak and off peak time. Scheduling problem is a well known optimization problem and the goal is to find the best possible optimal solution. In this paper, we used integer linear programming technique for scheduling problem in a high level of synthesis. In addition, time and resource constrained scheduling was used. An optimal result was fully obtained by using the software GLPK/AMPL. This model can be adopted to solve other scheduling problems, such as the Lecture Theatre, Cinemas and Work Shift.

## 1. INTRODUCTION

Hospital Operating Theatre (OT) scheduling involves an arrangement of several operating rooms to the medical surgeons in a period of time. In the health service sector such as government or private hospitals, the scheduling of Operating Theatre plays an important role towards achieving their goals. Their main goal is to meet the patients' satisfaction by minimizing his/her total waiting time before undergoing major or minor operations [3]. Poor scheduling of Operating Theatre may cause longer waiting time and can also worsen the patient's disease. In this case, an effective schedule has to be developed in order to improve the reputation and performance of government as well as private hospitals.

Another objective is to minimize the costs involved with Operating Theatre such as staffing cost, or also known as overtime cost, which is greater than that of regular working hours [3]. For example, a hospital needs to hire part-time medical specialists to perform major operations during weekends or outside office hours. As a result, the hospital has to spend more to satisfy the needs of the patient. This has led to the need of optimizing the availability of the operating theatre by minimizing the under-allocation of each operation room. Several recent researches have focused on both the objectives of minimizing the overtime cost and maximizing the utilization of operating theatre [1, 2, 3, 5, 9, 10].

Scheduling of Hospital Operating Theatre can be divided into several types of strategies which are open scheduling, block scheduling and modified block scheduling [4]. Open scheduling allows surgical cases to be assigned to an operating room available at the convenience of the surgeons [3, 6]. Meanwhile, for block scheduling, specific surgeons or groups of surgeons are assigned a set of time blocks, normally for some weeks or months, into which they can arrange their surgical cases. In the pure form, the surgeon or group "owns" their time blocks. None of those time blocks can be released [4]. "Modified block scheduling" is modified into two ways to increase its flexibility. Either some time is blocked and some is left open, or unused block time is released at an agreed-upon time before surgery, for instance 72 hours [4].

Implementation of different types of operating theatre scheduling is based on the complexities of the real situation in the hospitals. There are several techniques that have been used by past researchers in order to achieve their objectives and reduce the constraints such as linear programming [5], heuristic algorithm [6], hybrid genetic algorithm [4] and mixed integer linear programming [7].

In this study, a procedure to optimize the utilization of operating theatre by minimizing the total deviation from the weekly target hours and relative amount of under-allocation of rooms to different departments is pretested. Using a small data set, an integer linear programming technique was used to find the optimal operating theatre schedule in one-week period with open scheduling strategy.

The rest of the paper is organized as follows. Section 2 describes the numerical example and model formulation. Section 3 presents the result obtained to evaluate the performance of an integer linear programming technique. Section 4 concludes the paper.

## 2. NUMERICAL EXAMPLE AND MODEL FORMULATION

This study is a simplified version of the real operation scheduling problem. Consider a hospital that has 10-staff operating rooms serving 6 departments: surgery, gynaecology, ophthalmology, otolaryngology, oral surgery and emergency. There are 8 main surgical rooms and 2 elective outpatient surgery (EOPS) rooms. An operating room is either "short hours" or "long hours", depending upon the daily number of hours the room is in use. Because of the socialized nature of health care, all surgeries are scheduled during work days only from Monday through Friday.

Table 1. Surgery Room Availability.

| Availability hours Weekday | Main short (hours) | Main long (hours) | EOPS short (hours) | EOPS long (hours) |
|---|---|---|---|---|
| Monday | 7.5 | 9.0 | 7.5 | 8.0 |
| Tuesday | 7.5 | 9.0 | 7.5 | 8.0 |
| Wednesday | 7.5 | 9.0 | 7.5 | 8.0 |
| Thursday | 7.5 | 9.0 | 7.5 | 8.0 |
| Friday | 7.5 | 9.0 | 7.5 | 8.0 |
| **Number of rooms** | **4 per day** | **4 per day** | **1 per day** | **1 per day** |

Table 2. Weekly Demand for Operating Rooms Hours.

| Department | Weekly target hours ($h_j$) | Allowable Limits of Under-allocated hours ($u_j$) |
|---|---|---|
| Surgery | 187.0 | 10.0 |
| Gynaecology | 117.4 | 10.0 |
| Ophthalmology | 39.4 | 10.0 |
| Oral surgery | 19.9 | 10.0 |
| Otolaryngology | 26.3 | 10.0 |
| Emergency | 5.4 | 3.0 |

Table 1 summarizes the daily availability of the different types of rooms and Table 2 provides the weekly demand for operating room hours. The limit on the under-allocated hours in Table 2 is the most a department can be denied relative to its weekly request. We can devise a daily schedule that most satisfies the weekly target hours for different departments. We set target hours for each department as a goal. The given situation involves 6 departments and types of rooms.

Let say,

$x_{ijk}$ = number of rooms of type $i$ assigned to department $j$ on k days {k=1, 2,...5}

$u_j$ = maximum number of under-allocated hours allowed in department $j$

$h_j$ = requested ideal target hours for department $j$

$a_{ik}$ = number of rooms of type $i$ available on a day

$d_{ik}$ = duration of availability in hours of room type $i$ on day $k$

## 2.1 Fitness Function

The objective of the model is to minimize the total deviation from the weekly target hours and relative amount of under-allocation of rooms to different departments. The problem denotes that the ratio $\frac{s_j}{h_j}$, measures the relative amount of under-allocation for the department $j$.

$$Minimize \sum \frac{s_j}{h_j} \qquad (1)$$

where $s_j$ is the number of under-allocated hours in each department $j$.

## 2.2 Constraints

We need to assign each department with a constraint:

$$s_j \leq u_j \tag{2}$$

where the number of under-allocated hours has to be less or equal than the allowable limit. For each type of room $i$ assigned to department $j$ on each day $k$, we set the constraint as:

$$\sum x_{ijk} \leq a_{ik} \tag{3}$$

where the sum of the rooms that is assigned to each department on a particular day has to be less or equal to the number of rooms of type $i$ available on that day. We add up another restriction:

$$s_j \geq h_j - \sum d_{ik} x_{ijk} \tag{4}$$

in order to reach the weekly target hours. All parameters and variables are positive:

$$s_j, x_{ijk}, h_j, u_j, a_{ik}, d_{ik} \geq 0 \tag{5}$$

This present study involves small decision variables, constraints and parameters. The model searches for an optimal solution using an integer programming algorithm. The coding was programmed using GLPK/AMPL and the optimal solution for a 5-day schedule was obtained less than 2.0 seconds running on a 2.26GHz PC. This is very fast compared to other approaches, such as heuristic Genetic Algorithms, that have been used elsewhere [8]. The efficiency of the Linear Programming approach in solving the scheduling problem has been mentioned in many past studies [5].

## 3. NUMERICAL RESULTS

The results are presented in Table 3 (integer approach) and Table 4 (real number approach). Each table satisfies the goal and sets out to detail the allocation of rooms by type to different departments during the working week which is from Monday through Friday.

Table 3. Room Allocation in Integer Solution.

|  | Main short ||||| Main long ||||| EOPS short ||||| EOPS long |||||
|---|---|---|---|---|---|---|---|---|---|---|---|---|---|---|---|---|---|---|---|---|
|  | M | T | W | T | F | M | T | W | T | F | M | T | W | T | F | M | T | W | T | F |
| Surgery | 0 | 0 | 4 | 2 | 4 | 0 | 4 | 4 | 0 | 3 | 0 | 0 | 0 | 0 | 0 | 1 | 0 | 0 | 1 | 0 |
| Gynaecology | 4 | 4 | 0 | 1 | 0 | 1 | 0 | 0 | 2 | 0 | 1 | 0 | 1 | 0 | 0 | 0 | 0 | 1 | 0 | 0 |
| Ophthalmology | 0 | 0 | 0 | 1 | 0 | 0 | 0 | 0 | 2 | 1 | 0 | 0 | 0 | 0 | 0 | 0 | 1 | 0 | 0 | 0 |
| Oral surgery | 0 | 0 | 0 | 0 | 0 | 0 | 0 | 0 | 0 | 0 | 0 | 1 | 0 | 0 | 1 | 0 | 0 | 0 | 0 | 1 |
| Otolaryngology | 0 | 0 | 0 | 0 | 0 | 3 | 0 | 0 | 0 | 0 | 0 | 0 | 0 | 0 | 0 | 0 | 0 | 0 | 0 | 0 |
| Emergency | 0 | 0 | 0 | 0 | 0 | 0 | 0 | 0 | 0 | 0 | 0 | 0 | 0 | 1 | 0 | 0 | 0 | 0 | 0 | 0 |

Table 4. Room Allocation in Real Number Solution.

|  | Main short ||||| Main long ||||| EOPS short ||||| EOPS long |||||
|---|---|---|---|---|---|---|---|---|---|---|---|---|---|---|---|---|---|---|---|---|
|  | M | T | W | T | F | M | T | W | T | F | M | T | W | T | F | M | T | W | T | F |
| Surgery | 0 | 4 | 4 | 4 | 0 | 0 | 0 | 0 | 1.79 | 4 | 0 | 0 | 0 | 1 | 1 | 1 | 1 | 1 | 1 | 1 |
| Gynaecology | 0 | 0 | 0 | 0 | 4 | 0 | 1 | 4 | 2.21 | 0 | 1 | 1 | 1 | 0 | 0 | 0 | 0 | 0 | 0 | 0 |
| Ophthalmology | 4 | 0 | 0 | 0 | 0 | 0.48 | 0.79 | 0 | 0 | 0 | 0 | 0 | 0 | 0 | 0 | 0 | 0 | 0 | 0 | 0 |
| Oral surgery | 0 | 0 | 0 | 0 | 0 | 0 | 2.21 | 0 | 0 | 0 | 0 | 0 | 0 | 0 | 0 | 0 | 0 | 0 | 0 | 0 |
| Otolaryngology | 0 | 0 | 0 | 0 | 0 | 2.92 | 0 | 0 | 0 | 0 | 0 | 0 | 0 | 0 | 0 | 0 | 0 | 0 | 0 | 0 |
| Emergency | 0 | 0 | 0 | 0 | 0 | 0.6 | 0 | 0 | 0 | 0 | 0 | 0 | 0 | 0 | 0 | 0 | 0 | 0 | 0 | 0 |

From table 3 and 4, each room meets the requirement per day and fulfills the weekly target hours. Integer solution shows better results than real number solution because it produces a better and logical use of the available rooms. The percentages of departmental allocation by integer solution are: surgery 102%, gynaecology 100%, ophthalmology 108%, oral surgery 116%, otolaryngology 103% and emergency 139%. The logic of the model is that it may not be possible to satisfy the target hours for department *j* without applying the surplus. The objective is to determine a schedule that minimizes the relative amount of under-allocation of rooms to the different departments. The nonnegative variables *s-j* and *s+j* represent the under and over allocation of hours relative to the target for department *j*. The ratio $\frac{s_{-j}}{h_j}$ measures the relative amount of under-allocated slack to department *j*. Table 5, shows the actual value of surplus that is used in the problem and is compared to the allowable limits of under-allocated hours.

Table 5. Surplus and Allowable Limits of Under-allocated hours.

| Department | Surplus (hours) ($s_j$) | Percentages of departmental allocation (%) | Allowable Limits of Under-allocated hours ($u_j$) |
|---|---|---|---|
| Surgery | 3.0 | 102 | 10.0 |
| Gynaecology | 0.1 | 100 | 10.0 |
| Ophthalmology | 3.1 | 108 | 10.0 |
| Oral surgery | 3.1 | 116 | 10.0 |
| Otolaryngology | 0.7 | 103 | 10.0 |
| Emergency | 2.1 | 139 | 3.0 |

## 4. CONCLUSION

This problem is very common in managing operating rooms for government and private hospitals. Many studies have been done to find the best schedule which is able to meet the patient's satisfaction and minimize the total cost of the operating theatre. This paper uses the integer linear programming approach to solve the operating theatre problem with open scheduling strategy. This study is concerned more on optimal utilization of the availability of the operating theatre. Results show that the integer linear programming technique managed to obtain the optimal schedule that meets the objectives of the study and also the related constraints. Moreover, we can provide extra treatment on surgery, ophthalmology, oral surgery, otolaryngology and emergency if needed based on the availability of the rooms. The scheduling procedure is efficient where we manage to meet all the requirements and give a major impact on the productivity of the hospital treatment. It greatly outperforms older manual scheduling methods and optimizes real-time workloads. Further study can be done by focusing on both the objectives to maximize the utilization of the Operating Theatre as well as minimize the total cost with a large set of sample data that represents the real situation in the hospitals. Also, different approaches can be used, such as heuristic algorithm technique, which are able to solve NP-hard scheduling problem for other related areas.

**Acknowledgments**
Funding: University Tun Hussein Onn Malaysia, Parit Raja, Batu Pahat, Johor, Malaysia.